\documentclass{cmslatex}
\usepackage{graphicx,amsmath,amsfonts,amssymb,bbm}
\usepackage{psfrag}
\usepackage{epsf}
\usepackage{hyperref}

\usepackage[paperwidth=7.25in, paperheight=10in, margin=.875in]{geometry}

\def\theequation {\thesection.\arabic{equation}}
\makeatletter\@addtoreset {equation}{section}\makeatother

\newtheorem{theo}{Theorem}[section]
\newtheorem{lem}[theo]{Lemma}

\newtheorem{cor}[theo]{Corollary}
\newtheorem{rem}[theo]{Remark}

\newenvironment{Proof}
{\begin{trivlist} \item[]{\bf Proof. }}%
{\hspace*{\fill}$\rule{.3\baselineskip}{.35\baselineskip}$\end{trivlist}}

\renewcommand{\geq}{\geqslant}
\renewcommand{\leq}{\leqslant}

\renewcommand{\phi}{\varphi}
\newcommand{\be}{\begin{eqnarray}}
\newcommand{\ee}{\end{eqnarray}}

\begin{document}

\title{\bf On the linearized log--KdV equation}

\author{Dmitry E. Pelinovsky\thanks{Department of Mathematics and Statistics, McMaster University, Hamilton, Ontario, Canada, L8S 4K1; e-mail: dmpeli@math.mcmaster.ca}}

   \pagestyle{myheadings}
    \markboth{On the linearized log--KdV equation}{Dmitry E. Pelinovsky.}
\date{\today}
\maketitle

\begin{abstract}
The logarithmic KdV (log--KdV) equation admits global solutions in an energy space
and exhibits Gaussian solitary waves. Orbital stability of Gaussian solitary waves
is known to be an open problem. We address properties of solutions to the linearized
log--KdV equation at the Gaussian solitary waves. By using the decomposition
of solutions in the energy space in terms of Hermite functions, we show
that the time evolution is related to a Jacobi difference operator with a limit
circle at infinity. This exact reduction allows us to characterize both
spectral and linear orbital stability of solitary waves. We also introduce a convolution
representation of solutions to the log--KdV equation with the Gaussian weight and
show that the time evolution in such a weighted space is dissipative with the exponential rate of decay.
\end{abstract}

\section{Introduction}

We address the logarithmic Korteweg--de Vries (log-KdV) equation derived
in the context of solitary waves in granular chains with Hertzian interaction forces \cite{Chat,DP14,JP13}:
\begin{equation}
\label{logKdV}
v_t + v_{xxx} + (v \log|v|)_x = 0, \quad (x,t) \in \mathbb{R} \times \mathbb{R}.
\end{equation}
The log--KdV equation (\ref{logKdV}) has a two-parameter family of
Gaussian solitary waves
\begin{equation}
\label{soliton-orbit}
v(t,x) = e^{c} V(x-ct-a), \quad a,c \in \mathbb{R},
\end{equation}
where $V$ is a symmetric standing wave given by
\begin{equation}
\label{Gaussian}
V(x) := e^{\frac{1}{2} - \frac{x^2}{4}}, \quad x \in \mathbb{R}.
\end{equation}

Global solutions to the log--KdV equation (\ref{logKdV}) were constructed in \cite{CP14}
in the energy space
\begin{equation}
\label{energy-space}
X := \left\{ v \in H^1(\mathbb{R}) : \quad v^2 \log|v| \in L^1(\mathbb{R}) \right\},
\end{equation}
by a modification of analytic methods available for the log--NLS equation \cite{Caz1}
(also reviewed in Section 9.3 in \cite{Caz}). In the energy space $X$, the following quantities
for the momentum and energy,
\begin{equation}
\label{momentum}
Q(v) = \int_{\mathbb{R}} v^2 dx
\end{equation}
and
\begin{equation}
\label{energy}
E(v) = \int_{\mathbb{R}} \left[ \left( \partial_x v \right)^2 - v^2 \log|v| + \frac{1}{2} v^2 \right] dx
\end{equation}
are non-increasing functions of time $t$. Uniqueness, continuous dependence, and energy conservation
are established in \cite{CP14} under the additional condition $\partial_x \log|v| \in L^{\infty}(\mathbb{R} \times \mathbb{R})$,
which is not satisfied in the neighborhood of the family of Gaussian solitary waves given by (\ref{soliton-orbit}) and (\ref{Gaussian}).
As a result, orbital stability of the Gaussian solitary waves was not established for the log--KdV equation (\ref{logKdV}), in a sharp
contrast with that in the log--NLS equation established in \cite{Caz2}.

A possible path towards analysis of orbital stability of Gaussian solitary waves is to study their linear and spectral stability
by using the linearized log--KdV equation
\begin{equation}
\label{linlogKdV}
u_t = \partial_x L u,
\end{equation}
where $L : H^2(\mathbb{R}) \cap L^2_2(\mathbb{R}) \to L^2(\mathbb{R})$
is the Schr\"{o}dinger operator with a harmonic potential given by
the differential expression
\begin{equation}
\label{Schrodinger}
L = - \partial_x^2 + \frac{1}{4} (x^2 - 6).
\end{equation}
The linearized log--KdV equation (\ref{linlogKdV}) arises at the formal linearization
of the log--KdV equation (\ref{logKdV}) at the perturbation $u := v - V$. The Schr\"{o}dinger
operator $L$ is the Hessian operator of the second variation of $E(v)$ at $v = V$. Although $E(v)$ in (\ref{energy}) is
not a $C^2$ functional at $v = 0$, the second variation of $E(v)$ is well defined at $v = V$ by
\begin{equation}
\label{energy-second-variation}
E_c(u) = \langle L u, u \rangle_{L^2},
\end{equation}
which is formally conserved in the time evolution of (\ref{linlogKdV}).

With new estimates to be obtained for the linearized log--KdV equation (\ref{linlogKdV}), we may hope to
develop an ultimate solution of the outstanding problem on the orbital stability of the Gaussian solitary waves.
Indeed, if we set $v(t,x) = V(x) + w(t,x)$ for the solution to the log--KdV equation (\ref{logKdV}),
we obtain an equivalent evolution equation
\begin{equation}
\label{logKdV-w}
w_t = \partial_x L w - \partial_x N(w),
\end{equation}
where the linearized part coincides with (\ref{linlogKdV}) and
the nonlinear term $N(w)$ is given by
$$
N(w) = w \log\left( 1 + \frac{w}{V} \right) + V \left[ \log \left( 1 + \frac{w}{V} \right) - \frac{w}{V} \right].
$$
It is clear that the nonlinear term $N(w)$ does not behave uniformly in $x$ unless
$w$ decays at least as fast as $V$ in (\ref{Gaussian}). On the other hand,
if $w(t,x) = V(x) h(t,x)$, where $h$ is a bounded function in its variables, then
$N(w) = V n(h)$, where $n(h) = h \log(1 + h) + \log(1 + h) - h$ is analytic in $h$ for any $h \in (-1,1)$.
Therefore, obtaining new estimates for the linearized log--KdV equation (\ref{linlogKdV})
in a function space with Gaussian weights may be useful in the nonlinear
analysis of the log--KdV equation (\ref{logKdV-w}).

The spectrum of $L$ in $L^2(\mathbb{R})$ consists of equally spaced simple eigenvalues
$$
\sigma(L) = \{ -1,0,1,2,\ldots\},
$$
which include exactly one negative eigenvalue with the eigenvector $V$ (defined without normalization).
Therefore, $E(v)$ is not convex at $V$ in $X$. Nevertheless, $E_c(u)$ is positive in the constrained space
\begin{equation}
\label{constrained-space}
X_c := \left\{ u \in H^1(\mathbb{R})\cap L^2_1(\mathbb{R}) : \quad \langle V, u \rangle_{L^2} = 0 \right\},
\end{equation}
which corresponds to the fixed value $Q(v) = Q(V)$ in (\ref{momentum}) at the linearized approximation.

Several results were obtained for the linearized log--KdV equation (\ref{linlogKdV}). In \cite{JP13},
{\em linear orbital stability} of Gaussian solitary waves was obtained in the following sense:
for every $u(0) \in X_c$, there exists a unique global solution $u(t) \in X_c$
of the linearized log--KdV equation (\ref{linlogKdV}) which satisfies the following bound
\begin{equation}
\label{global-bound}
\| u(t) \|_{H^1 \cap L^2_1} \leq C \| u(0) \|_{H^1\cap L^2_1}, \quad t \in \mathbb{R},
\end{equation}
for some $t$-independent positive constant $C$. This result was obtained in \cite{JP13}
from the conservation of $E_c(u)$
in the time evolution of smooth solutions to the linearized log--KdV equation (\ref{linlogKdV}),
the symplectic decomposition of the solution $u(t) \in X_c$, $t \in \mathbb{R}$
into the translational part and the residual part,
\begin{equation}
\label{symplectic-decomposition}
u(t) = b(t) \partial_x V + y(t), \quad \langle V, y(t) \rangle_{L^2} = \langle \partial_x^{-1} V, y(t) \rangle_{L^2} = 0,
\end{equation}
and the coercivity of $E_c(y)$ in the squared $H^1(\mathbb{R}) \cap L^2_1(\mathbb{R})$ norm in the sense
\begin{equation}
\label{coercivity-assumed}
\| y \|_{H^1 \cap L^2_1}^2 \leq C E_c(y), \quad \langle V, y \rangle_{L^2} = \langle \partial_x^{-1} V, y \rangle_{L^2} = 0,
\end{equation}
for some positive constant $C$. The first two facts are rather standard in energy methods for linear PDEs, whereas
the last fact, that is, the inequality (\ref{coercivity-assumed}), should not be taken as granted.

In \cite{CP14}, the nonzero spectrum of the linear operator
\begin{equation}
\label{linearized-operator}
\partial_x L : H^3(\mathbb{R}) \cap H^1_2(\mathbb{R}) \cap L^2_1(\mathbb{R})  \cap \dot{H}^{-1}(\mathbb{R}) \to L^2(\mathbb{R})
\end{equation}
was studied by using the Fourier transform that maps the third-order differential operator in physical space
into a second-order differential operator in Fourier space. Indeed, the Fourier transform
$\hat{u}(k) := \mathcal{F}(u)(k) = \int_{\mathbb{R}} u(x) e^{-ikx} dx$ applied to the linearized log--KdV equation
(\ref{linlogKdV}) yields the time evolution in the form
\begin{equation}
\label{linlogKdV-Fourier}
i \hat{u}_t = k \hat{L} \hat{u},
\end{equation}
where $\hat{L} : H^2(\mathbb{R}) \cap L^2_2(\mathbb{R}) \to L^2(\mathbb{R})$ is the Fourier image of operator
$L : H^2(\mathbb{R}) \cap L^2_2(\mathbb{R}) \to L^2(\mathbb{R})$ given by
\begin{equation}
\label{Schrodinger-Fourier}
\hat{L} = - \frac{1}{4} \partial_k^2 + k^2 - \frac{3}{2}.
\end{equation}
By reducing the eigenvalue problem for $k \hat{L}$ to the symmetric Sturm--Liouville form, it was found
in \cite{CP14} that the spectrum of $\partial_x L$ in $L^2(\mathbb{R})$
is purely discrete and consists of a double zero eigenvalue and a symmetric sequence of simple
purely imaginary eigenvalues $\{ \pm i \omega_n \}_{n \in \mathbb{N}}$ such that
$$
0 < \omega_1 < \omega_2 < ..., \quad \omega_n \to \infty \quad \mbox{\rm as} \quad n \to \infty.
$$
The double zero eigenvalue corresponds to the Jordan block
\begin{equation}
\label{null-space}
\partial_x L \partial_x V = 0, \quad \partial_x L V = - \partial_x V,
\end{equation}
whereas the purely imaginary eigenvalues $\lambda = \pm i \omega_n$ correspond to the
eigenfunctions $u = u_{\pm n}(x)$, which are smooth in $x$ but decay
algebraically as $|x| \to \infty$. The Fourier transform of $u_{\pm n}$ is supported on
the half-line $\mathbb{R}^{\pm}$ and decays like a Gaussian function at infinity.
It follows from the spectrum of $\partial_x L$ in $L^2(\mathbb{R})$
that the Gaussian solitary waves are {\em spectral stable}.

The eigenfunctions of $\partial_x L$ were also used in \cite{CP14} for spectral decompositions
in the constrained space $X_c$ in order to provide an alternative proof of the {\em linear orbital stability}
of the Gaussian solitary waves. This alternative technique still relies on the conjecture of
the coercivity of $E_c(y)$ in the squared $H^1(\mathbb{R}) \cap L^2_1(\mathbb{R})$ norm, that is,
on the inequality (\ref{coercivity-assumed}).

Because of the algebraic decay of the eigenfunctions of $\partial_x L$, it is not clear
if a function of $x$ that decays like the Gaussian function as $|x| \to \infty$
can be represented as series of eigenfunctions. Numerical simulations were undertaken in \cite{CP14}
to illustrate that solutions to the linearized log--KdV equation (\ref{linlogKdV})
with Gaussian initial data did not spread out as the time variable evolves.
Nevertheless, the solutions exhibited visible radiation at the left slopes.

The present work is developed to obtain new estimates for the linearized log--KdV equation (\ref{linlogKdV}).
In the first part of this work, we rely on the basis of Hermite functions in $L^2$-based Sobolev spaces and analyze
the discrete operators that replace the differential operators. In the second part, we obtain
dissipative estimates on the evolution of the linearized log--KdV equation (\ref{linlogKdV})
by representing solutions in terms of a convolution with the Gaussian solitary wave $V$.

The paper is structured as follows. Section 2 sets up the basic formalism of the Hermite functions
and reports useful technical estimates. Section 3 is devoted to the proof of the coercivity
bound (\ref{coercivity-assumed}). As explained above, this
coercivity bound implies {\em linear orbital stability} of the Gaussian solitary wave in the constrained
space $X_c$ and it is assumed to be granted in \cite{CP14,JP13}. The proof of coercivity
relies on the decomposition of $y$ in terms of the Hermite functions.

Section 4 is devoted to the analysis of linear evolution expressed in terms of the Hermite
functions. It is shown that this evolution reduces to the self-adjoint Jacobi difference operator
with the limit circle behavior at infinity. As a result, a boundary condition is needed
at infinity in order to define the spectrum of the Jacobi operator and to obtain the
norm-preserving property of the associated semi-group. Both {\em linear orbital stability} and {\em spectral stability}
of Gaussian solitary waves (\ref{Gaussian}) is equivalently proven by using the Jacobi difference operator.

In Section 5, we give numerical approximations of eigenvalues and eigenvectors of the Jacobi difference
equation. We show numerically that there exist subtle differences between the representation of
eigenvectors of $\partial_x L$ in the physical space and the representation
of these eigenvectors by using decomposition in terms of the Hermite functions.

Section 6 reports weighted estimates for solutions to the linearized log--KdV equation (\ref{linlogKdV})
by using a convolution representation with the Gaussian weight. We show that the convolution representation
is invariant under the time evolution of the linearized log--KdV equation (\ref{linlogKdV}),
which is expressed by a dissipative operator on a half-line. The semi-group
of the fundamental solution in the $L^2(\mathbb{R}) \cap L^{\infty}(\mathbb{R})$ norm
decays to zero exponentially fast as time goes to infinity.

Section 7 concludes the paper with discussions of further prospects.

\vspace{0.25cm}

{\bf Notations:} We denote with $H^s(\mathbb{R})$ the Sobolev space
of $s$-times weakly differentiable functions on the real line whose derivatives up to order $s$ are in $L^2(\mathbb{R})$.
The norm $\| u \|_{H^s}$ for $u$ in the Sobolev space $H^s(\mathbb{R})$ is equivalent to the  norm
$\| (I - \partial_x^2)^{s/2} u \|_{L^2}$ in the Lebesgue space $L^2(\mathbb{R})$. We denote with
$L^2_s(\mathbb{R})$ the space of square integrable functions with the weight $<x>^s = (1 + x^2)^{s/2}$.
The set $\mathbb{N}_0$ consists of all non-negative integers, whereas the set $\mathbb{N}$ includes
only positive integers. The sequence space $\ell^2(\mathbb{N})$ includes squared summable sequences,
whereas $\ell_0(\mathbb{N})$ contains finite (compactly supported) sequences.

\vspace{0.25cm}

{\bf Acknowledgements.} The author thanks Gerald Teschl and Thierry Gallay for help on obtaining
results reported in Sections 4 and 6, respectively. The research of the author is supported by
the NSERC Discovery grant.

\newpage

\section{Preliminaries}

We recall definitions of the Hermite functions \cite[Chapter 22]{AS}:
\begin{equation}
\label{Hermite-function} \phi_n(z) = \frac{1}{\sqrt{2^n n!
\sqrt{\pi}}} H_n(z) e^{-\frac{z^2}{2}}, \qquad n \in \mathbb{N}_0,
\end{equation}
where $\{ H_n \}_{n \in \mathbb{N}_0}$ denote the set of Hermite polynomials, e.g.,
\begin{eqnarray*}
H_0 & = & 1, \\
H_1 & = & 2z, \\
H_2 & = & 4 z^2 - 2, \\
H_3 & = & 8 z^3 - 12 z.
\end{eqnarray*}
Hermite functions satisfy the Schr\"{o}dinger equation for a quantum harmonic oscillator:
\begin{equation}
\label{Schrodinger-Hermite}
-\phi_n''(z) + z^2 \phi_n(z) =  (1+2n) \phi_n(x), \quad n \in \mathbb{N}_0,
\end{equation}
at equally spaced energy levels. By the Sturm--Liouville theory \cite{Teschl},
the set of Hermite functions $\{ \phi_n \}_{n \in \mathbb{N}_0}$ forms an orthogonal
and normalized basis in $L^2(\mathbb{R})$.

In connection to the self-adjoint operator  $L : H^2(\mathbb{R}) \cap L^2_2(\mathbb{R}) \to L^2(\mathbb{R})$
given by (\ref{Schrodinger}), we obtain the eigenfunctions of $L u_n = (n-1) u_n$, $n \in \mathbb{N}_0$
from the correspondence $x = \sqrt{2} z$. With proper normalization, we define
\begin{equation}
\label{Hermite-function-ours} u_n(x) = \frac{1}{\sqrt{2^n n!
\sqrt{2 \pi}}} H_n\left(\frac{x}{\sqrt{2}}\right) e^{-\frac{x^2}{4}}, \qquad n \in \mathbb{N}_0.
\end{equation}
It follows from the well-known relations for Hermite polynomials
$$
H_n'(z) = 2n H_{n-1}(z), \quad 2 z H_n(z) = H_{n+1}(z) + 2 n H_{n-1}(z), \quad n \in \mathbb{N}_0,
$$
that functions in the sequence $\{ u_n \}_{n \in \mathbb{N}_0}$ satisfy the differential relations
\begin{equation}
\label{diff-Hermite}
2 u_n'(x) = -\sqrt{n+1} u_{n+1}(x) + \sqrt{n} u_{n-1}(x), \quad n \in \mathbb{N}_0.
\end{equation}

The following elementary result is needed in further estimates.

\vspace{0.25cm}

\begin{lem}
\label{lemma-f-sequence}
Let $\{ f_m \}_{m \in \mathbb{N}_0}$ be given by
$$
f_m = \prod_{k=1}^m \frac{\sqrt{2k-a}}{\sqrt{2k+b}}, \quad a,b \geq 0.
$$
Then, there is a positive constant $C$ such that
\begin{equation}
\label{bound-decay}
f_m \leq C m^{-(a+b)/4} \quad n \in \mathbb{N}.
\end{equation}
\end{lem}

\begin{Proof}
We write
\begin{equation}
\label{f-m-exp}
f_{m} = \exp\left[ \frac{1}{2} \sum_{k=1}^m \log\left(1 - \frac{a}{2k} \right) - \frac{1}{2} \sum_{k=1}^m \log\left(1 + \frac{b}{2k} \right)\right].
\end{equation}
By Taylor series, for every $k \in \mathbb{N}$ and every $a \in \mathbb{R}^+$, there is $C > 0$ such that
\begin{equation}
\label{Taylor}
\left| \log\left(1 - \frac{a}{2k} \right) +\frac{a}{2k} \right| \leq  \frac{C}{k^2}.
\end{equation}
Furthermore, we recall Euler's constant $\gamma \approx 0.577215$ given by the limit
\begin{equation}
\label{Euler-gamma}
\gamma := \lim_{m \to \infty} \left| \sum_{k=1}^{m-1} \frac{1}{k} - \log(m) \right|.
\end{equation}
Since $\sum_{k=1}^m k^{-2}$ is bounded as $m \to \infty$, the estimate (\ref{Taylor}) and
the limit (\ref{Euler-gamma}) yield the bound
\begin{equation}
\label{bound-log}
\left| \frac{1}{2} \sum_{k=1}^m \log\left(1 - \frac{a}{2k} \right) + \frac{a}{4} \log(m) \right| \leq C, \quad
\forall m \in \mathbb{N},
\end{equation}
for some positive constant $C$. Substituting (\ref{bound-log}) into (\ref{f-m-exp}) proves the desired bound (\ref{bound-decay}).
\end{Proof}

\vspace{0.25cm}

The following technical result is needed for the proof of coercivity of the energy function.

\vspace{0.25cm}

\begin{lem}
\label{lemma-projections}
Let $\{f_n\}_{n \in \mathbb{N}_0}$ by defined by $f_n := \langle \partial_x^{-1} u_0, u_n \rangle_{L^2}$.
Then, there is a positive constant $C$ such that
\begin{equation}
\label{projection-decay}
0 < f_n \leq C (1+n)^{-1/4}, \quad n \in \mathbb{N}_0.
\end{equation}
\end{lem}

\begin{Proof}
Multiplying the differential relation (\ref{diff-Hermite}) by $\partial_x^{-1} u_0$ and integrating by parts, we obtain
\begin{eqnarray}
\nonumber
\sqrt{n} \langle \partial_x^{-1} u_0, u_{n-1} \rangle_{L^2} - \sqrt{n+1} \langle \partial_x^{-1} u_0, u_{n+1} \rangle_{L^2}
& = & 2 \langle \partial_x^{-1} u_0, u_n' \rangle_{L^2} \\
& = & -2 \langle u_0, u_n \rangle_{L^2}, \quad n \in \mathbb{N}_0.
\label{difference-integral}
\end{eqnarray}
Integrating directly, we compute
\begin{equation}
\label{f-0}
f_0 =  \langle \partial_x^{-1} u_0, u_0 \rangle_{L^2} = \frac{1}{2} \| u_0 \|_{L^1}^2 = \sqrt{2\pi}.
\end{equation}
Furthermore, using (\ref{difference-integral}) at $n = 0$, we also compute
\begin{equation}
\label{f-1}
f_1 = \langle \partial_x^{-1} u_0, u_1 \rangle_{L^2} = 2 \| u_0 \|_{L^2}^2 = 2.
\end{equation}
Thanks to orthogonality of Hermite functions, the right-hand side of (\ref{difference-integral}) is zero for $n \in \mathbb{N}$
and the numerical sequence $\{ f_n \}_{n \in \mathbb{N}_0}$ satisfies the recurrence equation
\begin{equation}
\label{difference-second-order}
f_{n+1} = \frac{\sqrt{n}}{\sqrt{n+1}} f_{n-1}, \quad n \in \mathbb{N},
\end{equation}
starting with the initial values for $f_0$ and $f_1$ in (\ref{f-0}) and (\ref{f-1}).
The recurrence equation (\ref{difference-second-order}) admits the exact solution
\begin{equation}
\label{f-sequence}
f_{2m} = \left( \prod_{k=1}^m \frac{\sqrt{2k-1}}{\sqrt{2k}} \right) f_0, \quad
f_{2m+1} = \left( \prod_{k=1}^m \frac{\sqrt{2k}}{\sqrt{2k+1}} \right) f_1, \quad m \in \mathbb{N}.
\end{equation}
Applying the bound (\ref{bound-decay}) of Lemma \ref{lemma-f-sequence} with $a = 1$, $b = 0$ or $a = 0$, $b = 1$
yields the bound (\ref{projection-decay}).
\end{Proof}

\section{Coercivity of the energy function}

In order to prove the coercivity bound (\ref{coercivity-assumed}), we
define the $L$-compatible squared norm in space $H^1(\mathbb{R}) \cap L^2_1(\mathbb{R})$,
\begin{equation}
\label{norm-H-1}
\| u \|_{H^1 \cap L^2_1}^2 := \int_{\mathbb{R}} \left[ u_x^2 + \frac{1}{4} x^2 u^2 + \frac{1}{2} u^2 \right] dx
\end{equation}
The second variation $E_c(u)$ is defined by (\ref{energy-second-variation}).
The following theorem yields the coercivity bound for the energy function,
which was assumed in \cite{CP14,JP13} without a proof.

\vspace{0.25cm}

\begin{theo}
\label{theorem-coercivity}
There exists a constant $C \in (0,1)$ such that for every $y \in H^1(\mathbb{R}) \cap L^2_1(\mathbb{R})$
satisfying the constraints
\begin{equation}
\label{constraints}
\langle u_0, y \rangle_{L^2} = \langle \partial_x^{-1} u_0, y \rangle_{L^2} = 0,
\end{equation}
it is true that
\begin{equation}
\label{coercivity}
C \| y \|^2_{H^1 \cap L^2_1} \leq E_c(y) \leq \| y \|^2_{H^1 \cap L^2_1},
\end{equation}
where  $E_c(y) = \langle L y, y \rangle_{L^2}$.
\end{theo}

\vspace{0.25cm}

\begin{Proof}
The upper bound in (\ref{coercivity}) follows trivially from the identity
$$
E_c(y) + 2 \| y \|_{L^2}^2 = \| y \|_{H^1 \cap L^2_1}^2,
$$
whereas the lower bound holds if there is a constant $C > 0$ such that
for every $y \in  H^1(\mathbb{R}) \cap L^2_1(\mathbb{R})$ satisfying
constraints (\ref{constraints}), it is true that
\begin{equation}
\label{coercivity-L-2}
\| y \|_{L^2}^2 \leq C E_c(y).
\end{equation}
By the spectral theorem, we represent every $y \in H^1(\mathbb{R}) \cap L^2_1(\mathbb{R})$ by
\begin{equation}
\label{decomposition}
y = \sum_{n \in \mathbb{N}_0} c_n u_n, \quad c_n = \langle u_n, y \rangle_{L^2},
\end{equation}
where the vector $c := (c_0,c_1,c_2,...)$ belongs to $\ell^2_1(\mathbb{N}_0)$.
It follows from the first constraint in (\ref{constraints}) that $c_0 = 0$.
Using the norm in (\ref{norm-H-1}), we obtain
$$
E_c(y) = \sum_{n \in \mathbb{N}} (n-1) |c_n|^2 \geq \| y - c_1 u_1 \|_{L^2}^2.
$$
Therefore,
$$
\| y \|_{L^2}^2 = |c_1|^2 + \| y - c_1 u_1 \|_{L^2}^2 \leq |c_1|^2 + E_c(y),
$$
and coercivity (\ref{coercivity-L-2}) is proved if we can show that $|c_1|^2$ is bounded by $E_c(y)$ up to a multiplicative constant.
To show this, we use the second constraint in (\ref{constraints}).
Since $\langle \partial_x^{-1} u_0, u_1 \rangle_{L^2} = 2 \| u_0 \|_{L^2}^2 = 2$, as it follows from (\ref{f-1}),
we have
\begin{equation}
\label{projection-2}
2 c_1 = -\langle \partial_x^{-1} u_0, y - c_1 u_1 \rangle_{L^2} = -\sum_{n = 2}^{\infty} c_n \langle \partial_x^{-1} u_0, u_n \rangle_{L^2}.
\end{equation}
By Lemma \ref{lemma-projections}, there is a positive constant $C_0 > 0$ such that
\begin{equation}
\label{bound-C-0}
C_0 := \sum_{n =2}^{\infty} \frac{|\langle \partial_x^{-1} u_0, u_n \rangle_{L^2}|^2}{n-1} < \infty,
\end{equation}
which follows from convergence of $\sum_{n \in \mathbb{N}} n^{-3/2}$.
Hence, by using Cauchy--Schwarz inequality in (\ref{projection-2}), we obtain
$$
4 |c_1|^2 \leq C_0 \sum_{n = 2}^{\infty} (n-1) |c_n|^2 = C_0 E_c(y),
$$
so that the bound (\ref{coercivity-L-2}) follows. The statement of the theorem is proven.
\end{Proof}

\vspace{0.25cm}

\section{Time evolution of the linearized log--KdV equation}

The time evolution of the linearized log--KdV equation (\ref{linlogKdV}) is considered in the
constrained energy space $X_c$ given by (\ref{constrained-space}).
For a vector $c := (c_0,c_1,c_2,...) \in \ell^2_1(\mathbb{N}_0)$,
we use the decomposition involving the Hermite functions,
\begin{equation}
\label{decomposition-lin-log-KdV}
u(t) = \sum_{n \in \mathbb{N}_0} c_n(t) u_n, \quad c_n(t) = \langle u_n, u(t) \rangle_{L^2},
\end{equation}
By using $Lu_n = (n-1) u_n$ and the differential relations (\ref{diff-Hermite}), the evolution
problem for the vector $c \in \ell^2_1(\mathbb{N}_0)$ is written as the lattice differential equation
\begin{equation}
\label{lattice-c}
2 \frac{dc_n}{dt} = n \sqrt{n+1} c_{n+1} - (n-2) \sqrt{n} c_{n-1}, \quad n \in \mathbb{N}_0.
\end{equation}
It follows from (\ref{lattice-c}) for $n = 0$ that if $u(0) \in X_c$ (so that $c_0(0) = 0$),
then $c_0(t) = 0$ and $u(t) \in X_c$ for every $t$. If $c_0(t) = 0$, then it follows
from (\ref{lattice-c}) for $n = 1$ that the time evolution of a projection of $u(t)$ to $u_1$ (which is proportional to
the translational mode $\partial_x V$) is given by
\begin{equation}
\label{projection-c-1}
\frac{d c_1}{dt} = \frac{1}{\sqrt{2}} c_2.
\end{equation}
The projection $c_1(t)$ is decoupled from the rest of the system (\ref{lattice-c}).
Therefore, introducing $b_n = c_{n+1}$ for $n \in \mathbb{N}$, we
close the evolution system (\ref{lattice-c}) at the lattice differential equation
\begin{equation}
\label{lattice-b}
2 \frac{d b_n}{dt} = (n+1) \sqrt{n+2} b_{n+1} - (n-1) \sqrt{n+1} b_{n-1}, \quad n \in \mathbb{N}.
\end{equation}

Since $c \in \ell^2_1(\mathbb{N})$, then $b \in \ell^2_1(\mathbb{N})$, so that we can
introduce $a_n = \sqrt{n} b_n$, $n \in \mathbb{N}$ with the vector $a \in \ell^2(\mathbb{N})$.
The sequence $\{ a_n \}_{n \in \mathbb{N}}$ satisfies the evolution system in the skew-symmetric form
\begin{equation}
\label{lattice-a}
2 \frac{d a_n}{dt} = \sqrt{n(n+1)(n+2)} a_{n+1} - \sqrt{(n-1) n (n+1)} a_{n-1}, \quad n \in \mathbb{N}.
\end{equation}
The evolution system (\ref{lattice-a}) can be expressed in the symmetric form by using the transformation
\begin{equation}
\label{transformation-a}
a_n = i^n f_n, \quad n \in \mathbb{N}.
\end{equation}
The new sequence $\{f_n \}_{n \in \mathbb{N}}$ satisfies the evolution system
written in the operator form
\begin{equation}
\label{lattice-f}
\frac{d f}{d t} = \frac{i}{2} J f,
\end{equation}
where $J$ is the Jacobi operator defined by
\begin{equation}
\label{Jacobi-f}
(J f)_n := \sqrt{n(n+1)(n+2)} f_{n+1} + \sqrt{(n-1) n (n+1)} f_{n-1}, \quad n \in \mathbb{N}.
\end{equation}

The Jacobi difference equation is $J f = z f$. According to the definition in Section 2.6 of \cite{TeschlJacobi},
the Jacobi operator $J$ is said to have a limit circle at infinity if a solution $f$ of $Jf = zf$ with $f_1 = 1$
is in $\ell^2(\mathbb{N})$ for some $z \in \mathbb{C}$. By Lemma 2.15 in \cite{TeschlJacobi},
this property remains true for all $z \in \mathbb{C}$. The following lemma shows that
this is exactly our case.

\vspace{0.25cm}

\begin{lem}
\label{lemma-circle}
The Jacobi operator $J$ defined by (\ref{Jacobi-f}) has a limit circle at infinity.
\end{lem}

\vspace{0.25cm}

\begin{Proof}
Let us consider the case $z = 0$ and define a solution of $Jv = 0$ with $v_1 = 1$.
The numerical sequence $\{ v_n \}_{n \in \mathbb{N}}$ satisfies the recurrence relation
$$
v_{n+1} = -\frac{\sqrt{n-1}}{\sqrt{n+2}} v_{n-1}, \quad n \in \mathbb{N},
$$
starting with $v_1 = 1$. Then, $v_n = 0$ for even $n$, whereas
$v_n$ for odd $n$ is given by the exact solution
$$
v_{2m+1} = (-1)^m \prod_{k = 1}^m \frac{\sqrt{2k-1}}{\sqrt{2k+2}}, \quad m \in \mathbb{N}.
$$
By Lemma \ref{lemma-f-sequence} with $a = 1$ and $b = 2$, there exists a positive constant $C$ such that
\begin{equation}
\label{bound-v}
|v_{2m-1}| \leq C m^{-3/4} \quad m \in \mathbb{N}.
\end{equation}
This guarantees that $v \in \ell^2(\mathbb{N})$.
\end{Proof}

\vspace{0.25cm}

By Lemma 2.16 in \cite{TeschlJacobi}, the Jacobi operator $J_{\rm max} : D(J_{\rm max}) \to \ell^2(\mathbb{N})$ with the domain
\begin{equation}
\label{domain-J}
D(J_{\rm max}) := \{ f \in \ell^2(\mathbb{N}) : \quad J f \in \ell^2(\mathbb{N}) \}
\end{equation}
is self-adjoint if $W_{\infty}(f,g) = 0$ for all $f,g \in D(J_{\rm max})$, where
the discrete Wronskian is given by
\begin{equation}
\label{Wronskian}
W_n(f,g) := \sqrt{n(n+1)(n+2)} (f_n g_{n+1} - f_{n+1} g_n), \quad n \in \mathbb{N}
\end{equation}
and $W_{\infty}(f,g) = \lim_{n \to \infty} W_n(f,g)$.

In order to define a self-adjoint extension of the Jacobi operator $J$
with the limit circle at infinity, we need to define a boundary condition as follows:
\begin{equation}
\label{boundary-condition}
BC(J) := \left\{ v \in D(J_{\rm max}) : \quad W_{\infty}(v,f) = 0 \quad \mbox{\rm for some \;} f \in D(J_{\rm max}) \right\},
\end{equation}
where $v$ is real. By Lemma 2.17 and Theorem 2.18 in \cite{TeschlJacobi},
the operator $\mathcal{J} : \mathcal{D}(v) \to \ell^2(\mathbb{N})$,
where $v \in BC(J)$ and
\begin{equation}
\label{domain-condition}
\mathcal{D}(v) := \left\{ f \in D(J_{\rm max}) : \quad W_{\infty}(v,f) = 0 \right\},
\end{equation}
represents a self-adjoint extension of the Jacobi operator $J$ with the limit circle at infinity.

Moreover, by Lemma 2.19 in \cite{TeschlJacobi}, the real spectrum of $\mathcal{J}$ in $\ell^2(\mathbb{N})$
is purely discrete. Since there is at most one linearly independent solution of the Jacobi difference equation $Jf = zf$
thanks to the recurrence relation in (\ref{Jacobi-f}), each isolated eigenvalue of the real spectrum of $\mathcal{J}$ is simple.

By Lemma 2.20 in \cite{TeschlJacobi}, all self-adjoint extensions of $J_{\rm min} : \ell_0(\mathbb{N}) \to \ell^2(\mathbb{N})$
are uniquely defined by the choice $v \in BC(J)$ in (\ref{boundary-condition})
spanned by a linear combination of two linearly independent solutions of $Jf = 0$. Since the value of $f_0$ plays no role
for the Jacobi operator $J$ thanks again to the recurrence relation in (\ref{Jacobi-f}) and the value $f_1$ can be
uniquely normalized to $f_1 = 1$, we have a unique choice for $v$ given by the solution of $Jv = 0$ with $v_1 = 1$.
Combining these facts together, we have obtained the following result.

\vspace{0.25cm}

\begin{lem}
\label{lem-Jacobi}
Let $v \in BC(J)$ be a unique solution of $Jv = 0$ with $v_1 = 1$ determined in the proof of Lemma \ref{lemma-circle}.
Then, $\mathcal{J} : \mathcal{D}(v) \to \ell^2(\mathbb{N})$ with the domain in (\ref{domain-condition})
is a unique self-adjoint extension of the Jacobi operator $J$ given by (\ref{Jacobi-f}).
Moreover, the spectrum of $\mathcal{J}$ consists of a countable set of simple real isolated eigenvalues.
\end{lem}

\vspace{0.25cm}

The following theorem and corollary provide the {\em linear orbital stability} of Gaussian solitary waves (\ref{Gaussian})
expressed by using the decomposition in terms of Hermite functions. The same result was
obtained in \cite{CP14,JP13} by using alternative techniques involving either the energy method \cite{JP13}
or the spectral decompositions \cite{CP14}.

\vspace{0.25cm}

\begin{theo}
\label{theorem-linear-evolution}
For every $a(0) \in \ell^2(\mathbb{N})$, there exists a unique solution
$a(t) \in \ell^2(\mathbb{N})$ to the evolution system (\ref{lattice-a})
for every $t \in \mathbb{R}$ satisfying $\| a(t) \|_{\ell^2} = \| a(0) \|_{\ell^2}$.
\end{theo}

\vspace{0.25cm}

\begin{Proof}
The semi-group property of the solution operator $e^{\frac{i}{2} J t}$ both for $t \in \mathbb{R}^+$ and $t \in \mathbb{R}^-$
associated with the linear system (\ref{lattice-f})
follows from the result of Lemma \ref{lem-Jacobi} and the classical semi-group theory \cite{Pazy}.
The result is transferred to the sequence $a \in \ell^2(\mathbb{N})$ by using the transformation (\ref{transformation-a}).
\end{Proof}

\vspace{0.25cm}

\begin{cor}
For every $u(0) \in X_c$ given by (\ref{constrained-space}), there exists a unique solution
$u(t) \in X_c$ to the linearized log--KdV equation (\ref{linlogKdV})
for every $t \in \mathbb{R}^+$ satisfying $E_c(u(t)) = E_c(u(0))$.
\end{cor}

\vspace{0.25cm}

\begin{Proof}
By using the transformations $a_n = \sqrt{n} b_n$ and $b_n = c_{n+1}$ for $n \in \mathbb{N}$
and the decomposition (\ref{decomposition-lin-log-KdV}), we obtain
\begin{equation}
\label{equivalence-a-energy}
\| a \|_{\ell^2}^2 = \sum_{n \in \mathbb{N}} n |c_{n+1}|^2 = E_c(u),  \quad u \in X_c,
\end{equation}
where $c_0 = 0$ is set uniquely in $X_c$. Recall that $c_0 = 0$ is an invariant reduction
of the linearized log--KdV equation (\ref{linlogKdV}). The assertion of the corollary follows
from Theorem \ref{theorem-linear-evolution} and the equivalence (\ref{equivalence-a-energy}).
\end{Proof}

\vspace{0.25cm}

\begin{rem}
At the first glance, the detached equation (\ref{projection-c-1}) might imply that the projection $c_1(t)$
to the translational mode $u_1$ may grow at most linearly as $t \to \infty$ in the energy space $X_c$
with conserved $E_c(u(t)) = E_c(u(0))$. However, it follows directly from
the linearized log--KdV equation (\ref{linlogKdV}) \cite{CP14} for the solution $u \in C(\mathbb{R}^+,X_c)$ that
$$
\langle \partial_x^{-1} u_0, u(t) \rangle_{L^2} = \langle \partial_x^{-1} u_0, u(0) \rangle_{L^2}, \quad t \in \mathbb{R}.
$$
Therefore, we obtain as in the proof of Theorem \ref{theorem-coercivity} that
$$
2 c_1(t) = \langle \partial_x^{-1} u_0, u(0) \rangle_{L^2} - \sum_{n \in \mathbb{N}} c_{n+1}(t)
\langle \partial_x^{-1} u_0, u_n \rangle_{L^2},
$$
where both terms are globally bounded for all $t \in \mathbb{R}$.
\end{rem}

\vspace{0.25cm}

\section{Numerical approximations of eigenvalues and eigenvectors}

We discuss here numerical approximations of eigenvalues in the spectrum of
the self-adjoint operator $\mathcal{J} : \mathcal{D}(v) \to \ell^2(\mathbb{N})$ constructed
in Lemma \ref{lem-Jacobi}. Let $\{ z_k \}_{k \in \mathbb{Z}}$ denote the simple real eigenvalues
of $Jf = z f$ with the ordering
$$
\cdots < z_{-2} < z_{-1} < z_0 = 0 < z_1 < z_2 < \cdots
$$
These real eigenvalues of $\mathcal{J}$ transform to eigenvalues
$\lambda_k = \frac{i}{2} z_k \in i \mathbb{R}$ of the spectral problem $\partial_x L u = \lambda u$
by using the decomposition (\ref{decomposition-lin-log-KdV}). Therefore,
the result of Lemma \ref{lem-Jacobi} also provides an alternative proof of
the {\em spectral stability} of the Gaussian solitary waves (\ref{Gaussian}),
which is also established in \cite{CP14}. However, by comparing
the eigenvectors obtained in the two alternative approaches, we will see some
sharp differences in the definition of function spaces these eigenvectors belong to.

To proceed with numerical approximations, we note that if $v \in D(J_{\rm max})$ is a solution of
$Jv = 0$ like in the proof of Lemma \ref{lemma-circle}, then $v_n = 0$ for even $n$.
Let us denote $V_m = v_{2m-1}$ for $m \in \mathbb{N}$. From the bound (\ref{bound-v}),
we note that $V_m = \mathcal{O}(m^{-3/4})$ as $m \to \infty$.

\begin{figure}[htbp]
\center
\includegraphics[scale=0.4]{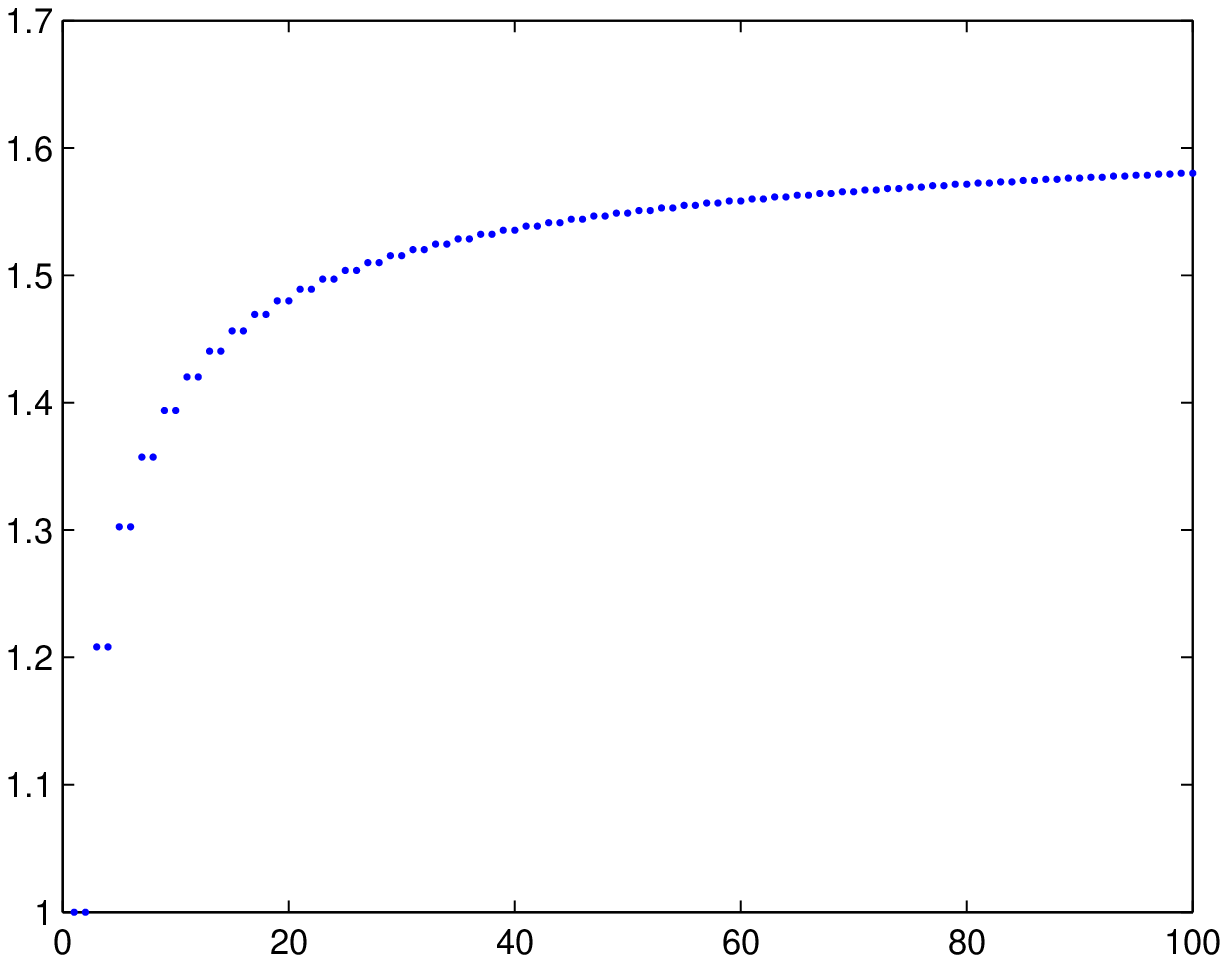}
\includegraphics[scale=0.4]{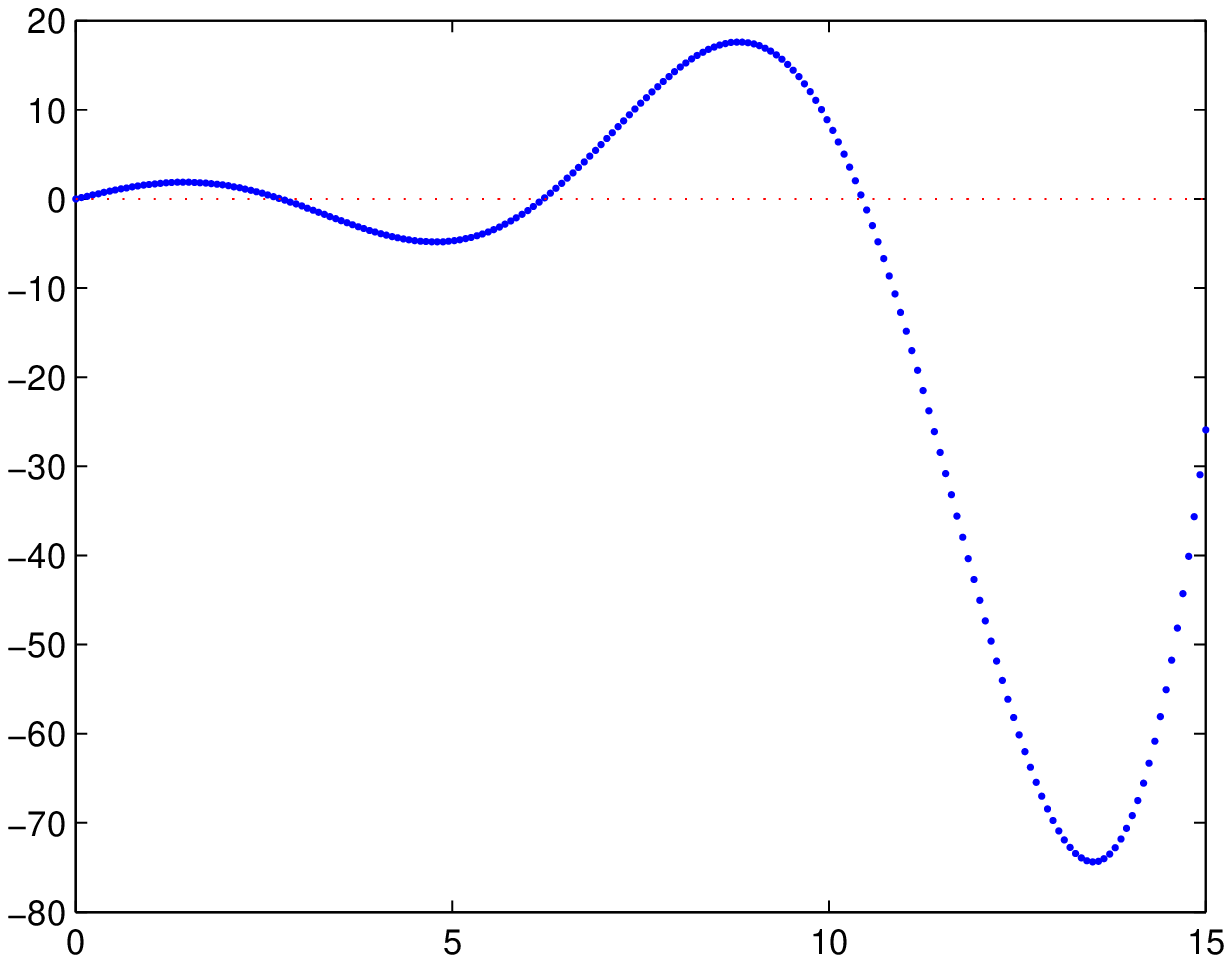}
\caption{(a) Convergence of the sequence $\{ W_n \}_{n \in \mathbb{N}}$ as $n \to \infty$ for $z = 1$. (b)
Oscillatory behavior of $W_{\infty}$ versus $z$.}
\label{fig-Wronskian}
\end{figure}

Let $f \in D(J_{\rm max})$ be a solution of $Jf = zf$ for $z \in \mathbb{R}^+$.
Then, we denote $A_m = f_{2m-1}$ and $B_m = f_{2m}$ for $m \in \mathbb{N}$.
It follows from the definition (\ref{Jacobi-f}) that $A$ and $B$ satisfy
the coupled system of difference equations:
\begin{eqnarray}
\label{coupled-difference}
\left\{ \begin{array}{l}
B_m = -\frac{\sqrt{2m-2}}{\sqrt{2m+1}} B_{m-1} + \frac{z}{\sqrt{(2m-1)(2m)(2m+1)}} A_m, \\
A_{m+1} = -\frac{\sqrt{2m-1}}{\sqrt{2m+2}} A_m + \frac{z}{\sqrt{(2m)(2m+1)(2m+2)}} B_m,
\end{array} \right. \quad m \in \mathbb{N},
\end{eqnarray}
starting with $A_1 = 1$. The discrete Wronskian (\ref{Wronskian}) is now explicitly computed as
\begin{equation}
\label{Wronskian-numerics}
W_n = \left\{ \begin{array}{l} \sqrt{(2m-1)(2m)(2m+1)} B_m V_m, \quad \quad \;\; n = 2m-1, \\
-\sqrt{(2m)(2m+1)(2m+2)} B_m V_{m+1}, \quad n = 2m \end{array} \right. \quad n \in \mathbb{N}.
\end{equation}
Since generally $B_m = \mathcal{O}(m^{-3/4})$ as $m \to \infty$, by applying Lemma \ref{lemma-f-sequence}
with $a = 2$ and $b = 1$ to the first equation of system (\ref{coupled-difference}),
the limit $W_{\infty} = \lim_{n \to \infty} |W_n|$ exists and is generally nonzero.
Moreover, the sign alternation of $\{ V_m \}_{m \in \mathbb{N}}$ and $\{ B_m \}_{m \in \mathbb{N}}$
ensures that the sequence $\{ W_n \}_{n \in \mathbb{N}}$ is sign-definite for large enough $n$,
so that the limit is actually $W_{\infty} = \lim_{n \to \infty} W_n$. This is confirmed in
Figure \ref{fig-Wronskian}(a), which shows the Wronskian sequence
$\{ W_n \}_{n \in \mathbb{N}}$ given by (\ref{Wronskian-numerics}) for $z = 1$.

Computing numerically $W_{\infty}$ by truncation of $\{ W_n \}_{n \in \mathbb{N}}$
at a sufficiently large $n$, e.g. at $n = 1000$, we plot
$W_{\infty}$ versus $z$ on Figure \ref{fig-Wronskian}(b). Oscillations of $W_{\infty}$ are observed
and the first two zeros of $W_{\infty}$ are located at
$$
z_1 \approx 2.7054, \quad z_2 \approx 6.1540.
$$
These values are nicely compared to the first two eigenvalues computed in \cite{CP14}
for $E_k = 2 z_k$:
$$
E_1 \approx 5.4109, \quad E_2 \approx 12.3080.
$$
The numerical approximations confirm that the eigenvalues obtained by using the Jacobi difference equation
are the same as the eigenvalues obtained in \cite{CP14} from the Sturm--Liouville problem derived in the Fourier space.
Numerically, we find for the first two zeros $z_{1,2}$ of the limiting Wronskian $W_{\infty}$ that
the decay rate of the sequence $\{ A_m \}_{m \in \mathbb{N}}$ remains generic:
$$
A_m = \mathcal{O}(m^{-3/4}) \quad \mbox{\rm as} \quad m \to \infty,
$$
but the decay rate of the sequence $\{ B_m \}_{m \in \mathbb{N}}$ becomes faster:
$$
B_m = \mathcal{O}(m^{-5/4})  \quad \mbox{\rm as} \quad m \to \infty.
$$

Let us now recall the correspondence of eigenvectors of the Jacobi difference equation
$Jf = zf$ and eigenvectors of the linearized log--KdV operator $\partial_x Lu = \lambda u$.
From the previous transformations, we obtain
\begin{equation}
\label{decomposition-numerics}
y := u - c_1 u_1 = \sum_{n \in \mathbb{N}} b_n u_{n+1} = \sum_{n \in \mathbb{N}} \frac{i^n}{\sqrt{n}} f_n u_{n+1}
= y_{\rm odd} + i y_{\rm even},
\end{equation}
where
\begin{equation}
\label{y-odd-even}
y_{\rm odd} = \sum_{m \in \mathbb{N}} \frac{(-1)^m}{\sqrt{2m}} B_m u_{2m+1} \quad {\rm and} \quad
y_{\rm even} = \sum_{m \in \mathbb{N}} \frac{(-1)^{m-1}}{\sqrt{2m-1}} A_m u_{2m}
\end{equation}
are respectively the odd and even components of the eigenvector with respect to $x$. Thanks to the decay of
the sequences $\{ A_m \}_{m \in \mathbb{N}}$ and $\{ B_m \}_{m \in \mathbb{N}}$, we note that
\begin{equation}
\label{y-spaces}
y_{\rm odd} \in H^2(\mathbb{R}) \cap L^2_2(\mathbb{R}), \quad y_{\rm even} \in H^1(\mathbb{R}) \cap L^2_1(\mathbb{R}),
\end{equation}
but that $y_{\rm even} \notin H^2(\mathbb{R}) \cap L^2_2(\mathbb{R})$. Therefore,
generally $y \notin {\rm Dom}(\partial_x L)$ defined by (\ref{linearized-operator}).
Thus, the eigenvector $y$ given by (\ref{decomposition-numerics})
does not solve the eigenvalue problem $\partial_x L y = \lambda y$ in the
classical sense compared to the eigenvectors constructed in \cite{CP14} with the Fourier transform.

In order to clarify the sense for the eigenvectors given by (\ref{decomposition-numerics}),
we denote $\lambda = \frac{i}{2} z$ and project the eigenvalue problem $\partial_x L u = \lambda u$
to $u_1$ and $y$. The projection $c_1$ is uniquely found by
\begin{equation}
\label{projection-eigenvalue-c-1}
z c_1 = \sqrt{2} A_1,
\end{equation}
which can also be obtained from the projection equation (\ref{projection-c-1}).
The component $y$ satisfies formally $\lambda y = \partial_x L y$. After separating
the even and odd parts of the eigenvalue problem, we obtain the coupled system
\begin{equation}
\label{linearized-formal}
z y_{\rm odd} = 2 \partial_x L y_{\rm even}, \quad -z y_{\rm even} = 2 \partial_x L y_{\rm odd}.
\end{equation}
As we have indicated above, it is difficult to prove that each term
of the coupled system (\ref{linearized-formal}) belongs to $L^2(\mathbb{R})$ if $y_{\rm even}$ and $y_{\rm odd}$
are given by the decomposition (\ref{decomposition-numerics}) in terms of the Hermite functions.
In order to formulate the coupled problem (\ref{linearized-formal}) rigorously, we would like to show that
the components of the eigenvector belong to
\begin{equation}
\label{function-space}
y_{\rm odd} \in H^2(\mathbb{R}) \cap L^2_2(\mathbb{R}) \cap H^{-1}(\mathbb{R}), \quad
y_{\rm even} \in H^1(\mathbb{R}) \cap L^2_1(\mathbb{R}) \cap H^{-1}(\mathbb{R}),
\end{equation}
and satisfy the coupled system
\begin{equation}
\label{linearized-classical}
z L^{-1} \partial_x^{-1} y_{\rm odd} = 2 y_{\rm even}, \quad -z \partial_x^{-1} y_{\rm even} = 2 L y_{\rm odd},
\end{equation}
where each term of system (\ref{linearized-classical}) is now defined in $L^2(\mathbb{R})$.

To show (\ref{function-space}) and (\ref{linearized-classical}), we proceed as follows.
According to (\ref{y-odd-even}) and (\ref{y-spaces}), $y_{\rm odd} \in H^2(\mathbb{R}) \cap L^2_2(\mathbb{R})$
is odd, hence $\partial_x^{-1} y_{\rm odd} \in L^2(\mathbb{R})$ and the first constraint in (\ref{function-space})
is satisfied. Since the kernel of $L$ is spanned by the odd function, we have $\partial_x^{-1} y_{\rm odd} \in {\rm Range}(L)$
so that $L^{-1} \partial_x^{-1} y_{\rm odd} \in L^2(\mathbb{R})$ and $y_{\rm even} \in L^2(\mathbb{R})$ as is given by the first
equation in (\ref{linearized-classical}). Similarly, from (\ref{y-odd-even}) and (\ref{y-spaces}),
we have $L y_{\rm odd} \in L^2(\mathbb{R})$ so that the second equation in (\ref{linearized-classical})
implies that $\partial_x^{-1} y_{\rm even} \in L^2(\mathbb{R})$. Hence, the second constraint
in (\ref{function-space}) is satisfied. Thus, the coupled system (\ref{linearized-classical}) is well defined for the eigenvector
$y = y_{\rm odd} + i y_{\rm even}$ of the eigenvalue problem $\partial_x L y = \lambda y$
defined in the function space (\ref{function-space}).

Note that the formulation (\ref{linearized-classical}) also settles the issue of zero eigenvalue $z_0 = 0$,
which should not be listed as an eigenvalue of the problem $\partial_z L y = \lambda y$.
Indeed, the first equation (\ref{linearized-classical}) with $z_0 = 0$ implies
$y_{\rm even} = 0$, hence $V_m = v_{2m-1} = 0$, where $v$ is a solution of $Jv = 0$.
Thus, the existence of the eigenvector $v \in \ell^2(\mathbb{N})$ for the eigenvalue $z_0 = 0$
of the Jacobi difference equation does not imply the existence
of the zero eigenvalue $\lambda = 0$ in the proper formulation (\ref{function-space})--(\ref{linearized-classical})
of the system $\partial_x L y = \lambda y$. The same result can be obtained from the projection
equation (\ref{projection-eigenvalue-c-1}). If $z_0 = 0$, then $V_1 \equiv A_1 = 0$, which
corresponds to the zero solution $v = 0$ of the Jacobi difference equation $Jv = 0$.

Delicate analytical issues in the decomposition (\ref{decomposition-numerics}) involving
Hermite functions are likely to be related to the fact that eigenvectors $u$ of the eigenvalue problem
$\partial_x L u = \lambda u $ decay algebraically as $|x| \to \infty$,
while the decay of each Hermite function $u_n$ in the decomposition (\ref{y-odd-even})
is given by a Gaussian function.

\vspace{0.25cm}

\section{Dissipative properties of the linearized log--KdV equation}

For the KdV equation with exponentially decaying solitary waves, the exponentially weighted spaces
were used to introduce effective dissipation in the long-time behavior of perturbations to the solitary waves
and to prove their asymptotic stability \cite{pw}. For the log--KdV equation with Gaussian solitary waves,
it makes sense to introduce Gaussian weights in order to obtain a dissipative evolution of the linear perturbations.
Here we show how the Gaussian weights can be introduced for the linearized log--KdV equation (\ref{linlogKdV}).

Let us represent a solution to the linearized log--KdV equation (\ref{linlogKdV}) in the following form
\begin{equation}
\label{convolution}
u(t,x) = a(t) u_0(x) + b(t) u_1(x) + y(t,x),
\end{equation}
with
\begin{equation}
\label{convolution-y}
y(t,x) = \int_{x}^{\infty} u_0(x') w(t,x-x') dx = \int_{-\infty}^0 u_0(x-z) w(t,z) dz,
\end{equation}
where $(a,b,w)$ are new variables to be found. It is clear that the representation
(\ref{convolution-y}) imposes restrictions on the class of functions of $y$ in
the energy space $X_c$. We will show that these restrictions are invariant
with respect to the time evolution of the linearized log--KdV equation (\ref{linlogKdV}).

We assume sufficient smoothness and decay of the variable $w$.
By using the explicit computation with $Lu_0 = -u_0$ and $x u_0(x) = -2 u_0'(x)$, we obtain
\begin{eqnarray*}
Ly & = & \int_{-\infty}^0 w(t,z) \left[ (Lu_0)(x-z) + \frac{1}{2} z (x-z) u_0(x-z) + \frac{1}{4} z^2 u_0(x-z) \right] dz \\
& = & \int_{-\infty}^0 w(t,z) \left[ -u_0(x-z) - z u_0'(x-z) + \frac{1}{4} z^2 u_0(x-z) \right] dz.
\end{eqnarray*}
Integrating by parts, we further obtain
\begin{eqnarray*}
\partial_x Ly & = & 2 w(t,0) u_0(x) + \int_{-\infty}^0 u_0(x-z) \left[ -w_z - (z w)_{zz} + \frac{1}{4} (z^2 w)_z \right] dz.
\end{eqnarray*}
Also recall that $\partial_x L u_0 = -u_0'(x) = \frac{1}{2} u_1(x)$ and $\partial_x L u_1 = 0$.
Bringing together the left-side and the right-side of the linearized log--KdV equation (\ref{linlogKdV})
under the decomposition (\ref{convolution})--(\ref{convolution-y}), we obtain the system of modulation equations
\begin{equation}
\label{modulation-equation}
\left\{ \begin{array}{l} \dot{a}(t) = 2 w(t,0), \\
\dot{b}(t) = \frac{1}{2} a(t), \end{array} \right.
\end{equation}
and the evolution problem
\begin{equation}
\label{heat}
w_t = H w,
\end{equation}
where the linear operator $H : {\rm Dom}(H) \to L^2(\mathbb{R}^-)$ with ${\rm Dom}(H) = \{ w \in L^2(\mathbb{R}^-), \quad Hw \in L^2(\mathbb{R}^-) \}$
is given by
\begin{equation}
\label{dissipative-operator}
(Hw)(z) = -z w_{zz} - 3 w_z + \frac{1}{4} (z^2 w)_z, \quad z < 0.
\end{equation}
Since $z = 0$ is a regular singular point of the differential operator $H : {\rm Dom}(H) \to L^2(\mathbb{R}^-)$,
no boundary condition is needed to be set at $z = 0$. We show that the differential
operator $H$ is dissipative in $L^2(\mathbb{R}^-)$.

\vspace{0.25cm}

\begin{lem}
\label{lemma-dissipative}
For every $w \in {\rm Dom}(H) \subset L^2(\mathbb{R}^-)$, we have
\begin{equation}
\label{H-form}
\langle Hw, w \rangle_{L^2(\mathbb{R}^-)} = - [w(0)]^2 + \int_{-\infty}^0 \left[ z (\partial_z w)^2 + \frac{1}{4} z w^2 \right] dz
\leq -\frac{1}{2} \| w \|_{L^2(\mathbb{R}^-)}^2.
\end{equation}
\end{lem}

\vspace{0.25cm}

\begin{Proof}
The proof of the equality follows from integration by parts for every $w \in {\rm Dom}(H)$:
{\small \begin{eqnarray*}
\int_{-\infty}^0 w \left[-z w_{zz} - 3 w_z + \frac{1}{4} (z^2 w)_z \right] dz
& = & \left[ -z w w_z - w^2 + \frac{1}{8} z^2 w^2 \right] \biggr|_{z \to -\infty}^{z = 0}
+ \int_{-\infty}^0 \left[ z (\partial_z w)^2 + \frac{1}{4} z w^2 \right] dz \\
& = & - [w(0)]^2
+ \int_{-\infty}^0 \left[ z (\partial_z w)^2 + \frac{1}{4} z w^2 \right] dz.
\end{eqnarray*}}
This yields the equality in (\ref{H-form}). The inequality in (\ref{H-form}) is proved from
the Younge inequality
$$
\| w \|_{L^2(\mathbb{R}^-)}^2 = -\int_{-\infty}^0 2 z w w_z dz \leq \alpha^2 \int_{-\infty}^0 |z| w_z^2 dz
+ \alpha^{-2} \int_{-\infty}^0 |z| w^2 dz,
$$
where $\alpha > 0$ is at our disposal. Picking $\alpha^2 = 2$ yields the inequality in (\ref{H-form}).
\end{Proof}

\vspace{0.25cm}

The semi-group theory for dissipative operators is fairly standard \cite{Pazy}, so we assume existence of
a strong solution to the evolution problem (\ref{heat}) for every $t > 0$. The next result shows
that this solution decays exponentially fast in the $L^2(\mathbb{R}^-)$ norm.

\vspace{0.25cm}

\begin{cor}
\label{cor-dissipation}
Let $w \in C(\mathbb{R}^+,{\rm Dom}(H)) \cap C^1(\mathbb{R}^+,L^2(\mathbb{R}^-))$ be a solution
of the evolution problem (\ref{heat}). Then, the solution satisfies
\begin{equation}
\label{decay-L-2}
\| w(t) \|_{L^2(\mathbb{R}^-)}^2 \leq \| w(0) \|_{L^2(\mathbb{R}^-)}^2 e^{-t}.
\end{equation}
\end{cor}

\vspace{0.25cm}

\begin{Proof}
The decay behavior (\ref{decay-L-2}) is obtained from a priori energy estimates. Indeed,
it follows from (\ref{H-form}) that
$$
\frac{1}{2} \frac{d}{dt} \| w(t) \|^2_{L^2(\mathbb{R}^-)} = \langle Hw, w \rangle_{L^2(\mathbb{R}^-)}
\leq -\frac{1}{2} \| w(t) \|^2_{L^2(\mathbb{R}^-)}.
$$
Gronwall's inequality yields the bound (\ref{decay-L-2}).
\end{Proof}

\vspace{0.25cm}

We recall that the solution $u \in X_c$ needs to satisfy the constraint $\langle u_0, u \rangle_{L^2} = 0$.
The constraint is invariant with respect to the time evolution of the linearized log--KdV equation (\ref{linlogKdV}).
These properties are equivalently represented in the decomposition (\ref{convolution})--(\ref{convolution-y}),
according to the following lemma.

\vspace{0.25cm}

\begin{lem}
\label{lemma-constraint}
For every $w \in {\rm Dom}(H) \subset L^2(\mathbb{R}^-)$, we have
\begin{equation}
\label{conserved-a}
a(t) + \int_{-\infty}^0 e^{-\frac{1}{8} z^2} w(t,z) dz = A, \quad t \in \mathbb{R}^+,
\end{equation}
where $A$ is constant in $t$. Moreover, if $u(0) \in X_c$, then $A = 0$ and
\begin{equation}
\label{decay-a}
|a(t)|^2 \leq \sqrt{\pi} \| w(0) \|_{L^2(\mathbb{R}^-)}^2 e^{-t}.
\end{equation}
\end{lem}

\vspace{0.25cm}

\begin{Proof}
We compute directly
\begin{eqnarray*}
\langle u_0, u(t) \rangle_{L^2} & = & a(t) + \frac{1}{\sqrt{2\pi}} \int_{-\infty}^{\infty} e^{-\frac{1}{4} x^2}
\left( \int_{-\infty}^0 e^{-\frac{1}{4} (x-z)^2} w(t,z) dz \right) dx \\
& = & a(t) + \frac{1}{\sqrt{2\pi}} \int_{-\infty}^{0} e^{-\frac{1}{8} z^2} w(t,z)
\left( \int_{-\infty}^{\infty} e^{-\frac{1}{2} \left( x - \frac{1}{2} z \right)^2} dx \right) dz \\
& = & a(t) + \int_{-\infty}^{0} e^{-\frac{1}{8} z^2} w(t,z) dz.
\end{eqnarray*}
Furthermore, $\langle u_0, u(t) \rangle_{L^2} = A$ is constant in $t$. As an alternative derivation,
one can compute from the evolution problem (\ref{heat}) that
$$
\frac{d}{dt} \int_{-\infty}^0 e^{-\frac{1}{8} z^2} w(t,z) dz = - 2 w(t,0)
$$
and then use the first modulation equation in system (\ref{modulation-equation}) for integration in $t$.

If $u(0) \in X_c$, then $u(t) \in X_c$ and $A = 0$. This yields $a(t)$ uniquely by
$$
a(t) = - \int_{-\infty}^0 e^{-\frac{1}{8} z^2} w(t,z) dz.
$$
Applying the decay bound (\ref{decay-L-2}) and the Cauchy--Schwarz inequality, we obtain
the decay bound (\ref{decay-a}).
\end{Proof}

\vspace{0.25cm}

\begin{cor}
\label{cor-modulation}
If $u(0) \in X_c$, then there is $b_{\infty} \in \mathbb{R}$
such that $b(t) \to b_{\infty}$ as $t \to \infty$.
\end{cor}

\vspace{0.25cm}

\begin{Proof}
It follows from the bound (\ref{decay-a}) that if $u(0) \in X_c$, then
$a(t) \to 0$ as $t \to \infty$. Since $a(t)$ decays to zero
exponentially fast, the assertion of the corollary follows
from integration of the second modulation equation in system (\ref{modulation-equation}).
\end{Proof}

\vspace{0.25cm}

Besides scattering to zero in the $L^2(\mathbb{R}^-)$ norm, the global solution
of the evolution problem (\ref{heat}) also scatters to zero in the $L^{\infty}(\mathbb{R}^-)$ norm.
The following lemma gives the relevant result based on a priori energy estimates.

\vspace{0.25cm}

\begin{lem}
\label{lemma-scattering}
Let $w$ be a smooth solution of the evolution problem (\ref{heat}) in a subset of $H^1(\mathbb{R}^-)$.
Then, there exist positive constants $\alpha$ and $C$ such that
\begin{equation}
\label{decay-L-infty}
\| w(t) \|_{L^{\infty}(\mathbb{R}^-)} \leq C \| w(0) \|_{H^1(\mathbb{R}^-)} e^{-\alpha t}.
\end{equation}
\end{lem}

\vspace{0.25cm}

\begin{Proof}
The proof is developed similarly to the estimates in Lemma \ref{lemma-dissipative} and Corollary \ref{cor-dissipation}
but the estimates are extended for $\| \partial_z w(t) \|_{L^2(\mathbb{R}^-)}$.
Differentiating (\ref{dissipative-operator}) in $z$, multiplying by $w_z$,
and integrating by parts, we obtain for smooth solution $w$:
{\small\begin{eqnarray*}
\langle w_z, (Hw)_z \rangle_{L^2(\mathbb{R}^-)} & = & \left[ -z w_z w_{zz} - \frac{3}{2} w_z^2
+ \frac{1}{8} z^2 w_z^2 + \frac{1}{4} w^2 \right] \biggr|_{z \to -\infty}^{z = 0}
+ \int_{-\infty}^0 z w_{zz}^2 dz + \frac{3}{4} \int_{-\infty}^0 z w_z^2 dz \\
& = & - \frac{3}{2} [ \partial_z w(0)]^2 + \frac{1}{4} [w(0)]^2
+ \int_{-\infty}^0 z w_{zz}^2 dz + \frac{3}{4} \int_{-\infty}^0 z w_z^2 dz.
\end{eqnarray*}}
As a result, smooth solutions to the evolution problem (\ref{heat}) satisfy
the differential inequality
\begin{eqnarray*}
\frac{d}{dt} \frac{1}{2} \| w_z\|_{L^2(\mathbb{R}^-)}^2 = \langle w_z, (Hw)_z \rangle_{L^2(\mathbb{R}^-)}
\leq \frac{1}{4} \left[ w(t,0) \right]^2 + \int_{-\infty}^0 z w_{zz}^2 dz + \frac{3}{4} \int_{-\infty}^0 z w_z^2 dz.
\end{eqnarray*}
By using Young's inequality, we estimate
$$
\left[ w(t,0) \right]^2 = 2 \int_{-\infty}^0 z w w_z dz \leq \beta^2 \| w_z\|_{L^2(\mathbb{R}^-)}^2 + \beta^{-2}
\| w\|_{L^2(\mathbb{R}^-)}^2
$$
and
$$
\| w_z\|_{L^2(\mathbb{R}^-)}^2 = -2 \int_{-\infty}^0 z w_z w_{zz} dz \leq \alpha^2 \int_{-\infty}^0 |z| w_{zz}^2 dz
+ \alpha^{-2} \int_{-\infty}^0 |z| w_z^2 dz,
$$
where $\alpha$ and $\beta$ are to our disposal. Picking $\alpha^2 = 2$ and assuming $\beta^2 < 2$,
we close the differential inequality as follows
\begin{eqnarray*}
\frac{d}{dt} \| w_z\|_{L^2(\mathbb{R}^-)}^2 \leq -\left(1 - \frac{\beta^2}{2}\right)\| w_z\|_{L^2(\mathbb{R}^-)}^2
+ \frac{1}{2 \beta^2} \| w\|_{L^2(\mathbb{R}^-)}^2.
\end{eqnarray*}
Thanks to the exponential decay in the bound (\ref{decay-L-2}), we can rewrite the differential inequality in the form
\begin{eqnarray*}
\frac{d}{dt} \left[ \| w_z\|_{L^2(\mathbb{R}^-)}^2 e^{\left(1 - \frac{\beta^2}{2}\right) t} \right] \leq
\frac{1}{2 \beta^2} \| w(0) \|_{L^2(\mathbb{R}^-)}^2 e^{-\frac{\beta^2}{2} t}.
\end{eqnarray*}
Integrating over time, we finally obtain
\begin{eqnarray*}
\| \partial_z w(t)\|_{L^2(\mathbb{R}^-)}^2 \leq \left( \| \partial_z w(0) \|_{L^2(\mathbb{R}^-)}^2
+ \beta^{-4} \| w(0) \|_{L^2(\mathbb{R}^-)}^2 \right) e^{-\left(1 - \frac{\beta^2}{2}\right) t},
\end{eqnarray*}
where $\beta^2 < 2$ is fixed arbitrarily. Thus, the $H^1(\mathbb{R}^-)$ norm of the smooth solution to the evolution problem (\ref{heat})
decays to zero exponentially fast as $t \to \infty$. The bound (\ref{decay-L-infty}) follows
by the Sobolev embedding of $H^1(\mathbb{R}^-)$ to $L^{\infty}(\mathbb{R}^-)$.
\end{Proof}

\vspace{0.25cm}

Combining the results of this section, we summarize the main
result on the dissipative properties of the solutions to the linearized log--KdV equation
(\ref{linlogKdV}) represented in the convolution form (\ref{convolution})--(\ref{convolution-y}).

\vspace{0.25cm}

\begin{theo}
\label{theorem-convolution}
Assume that the initial data $u(0) \in X_c$ is represented by the convolution form (\ref{convolution})--(\ref{convolution-y})
with some $a(0)$, $b(0)$, and $w(0) \in {\rm Dom}(H)$. There exists a solution of
the linearized log--KdV equation (\ref{linlogKdV}) represented in the convolution
form (\ref{convolution})--(\ref{convolution-y}) with unique $(a,b) \in C^1(\mathbb{R}^+,\mathbb{R}^2)$
and $w \in C(\mathbb{R}^+,{\rm Dom}(H)) \cap C^1(\mathbb{R}^+,L^2(\mathbb{R}^-))$.
Moreover, there is a $b_{\infty} \in \mathbb{R}$ such that
\begin{equation}
\label{scattering-result}
\lim_{t \to \infty} \| u(t) - b_0 u_1 \|_{L^2 \cap L^{\infty}}  = 0.
\end{equation}
\end{theo}

\vspace{0.25cm}

\begin{Proof}
The existence result follows from the existence of the semi-group to the evolution problem
(\ref{heat}) and the ODE theory for the system of modulation equations (\ref{modulation-equation}).
Since $u_0 \in L^1(\mathbb{R})$, the scattering result (\ref{scattering-result})
follows from the generalized Younge inequality, as well as the results of
Corollary \ref{cor-dissipation}, Lemma \ref{lemma-constraint},
Corollary \ref{cor-modulation}, and Lemma \ref{lemma-scattering}.
\end{Proof}

\vspace{0.25cm}

\section{Conclusion}

We have obtained new results for the linearized log--KdV equation.
By using Hermite function decompositions in Section 4, we have shown
analytically how the semi-group properties of the linear
evolution in the energy space can be recovered with the Jacobi difference operator.
We have also established numerically
in Section 5 the equivalence between computing the spectrum of the linearized operator with the Jacobi difference equation
and that with the differential equation.
Finally, we have used in Section 6 the convolution representation with the Gaussian weight to show that
the solution to the linearized log--KdV equation can decay to zero in the $L^2 \cap L^{\infty}$ norms.

It may be interesting to compare these results with the Fourier transform method used in the previous work \cite{CP14}.
From analysis of eigenfunctions of the spectral problem $\partial_x L u = \lambda u$, it is known that the eigenfunctions
are supported on a half-line in the Fourier space. The decomposition (\ref{decomposition-lin-log-KdV})
in terms of the Hermite functions in the physical space can be written equivalently as the decomposition
in terms of the Hermite functions in the Fourier space. The Jacobi difference equation
representing the spectral problem does not imply generally that the decomposition in the Fourier space returns an eigenfunction
supported on a half-line. This property is not explicitly seen in the computation of eigenvectors
with the Jacobi difference operator.

Another interesting observation is as follows. The linear evolution of the linearized log--KdV equation in the Fourier space (\ref{linlogKdV-Fourier})
can be analyzed separately for $k \in \mathbb{R}^+$ and $k \in \mathbb{R}^-$. 
Since the time evolution is given by the linear Schr\"{o}dinger-type equation,
the fundamental solution is norm-preserving in the energy space. If the Gaussian weight is introduced
on the positive half-line as follows:
$$
\hat{u}(t,k) = e^{-k^2} \hat{w}(t,k),\quad k > 0,
$$
then the time evolution is defined in the Fourier space by $i \hat{u}_t = \hat{H} \hat{u}$,
where the linear operator $\hat{H} : {\rm Dom}(\hat{H}) \to L^2(\mathbb{R}^+)$ is given by
\begin{equation}
\label{H-operator-Fourier}
\hat{H} = \frac{1}{4} k \partial_k^2 - k^2 \partial_k + k.
\end{equation}
If $H$ and $\hat{H}$ in (\ref{H-form}) and (\ref{H-operator-Fourier})
are extended on the entire line, then $H$ and $\hat{H}$ are Fourier images of each other.
Thus, a very similar introduction of the Gaussian weights (except, of course, the domains
in the physical and Fourier space) may result in either dissipative or norm-preserving solutions
of the linearized log--KdV equation.

Although the results obtained in this work give new estimates and new tools for analysis of the
linearized log--KdV equation, it is unclear in the present time how to deal with the main problem
of proving orbital stability of the Gaussian solitary waves in the nonlinear log--KdV equation.
This challenging problem will remain open to new researchers.

\end{document}